\documentclass[10pt]{amsart}
\usepackage{amsthm,amsmath,amssymb,amscd,graphics,enumerate, stmaryrd,xspace,verbatim, epic, eepic,url}

\usepackage[T1]{fontenc}
\usepackage[usenames,dvipsnames]{color}

\usepackage[colorlinks=true]{hyperref}

\usepackage[active]{srcltx}
\usepackage[all]{xypic}
\SelectTips{cm}{}

{
   \newtheorem{theorem}[subsubsection]{Theorem}
      \newtheorem*{theorem*}{Theorem}
   \newtheorem{proposition}[subsubsection]{Proposition}

   \newtheorem{corollary}[subsubsection]{Corollary}
   
   \newtheorem*{conjecture*}{Conjecture}

}
{\theoremstyle{definition}
          \newtheorem*{exercise*}{Exercise}

   \newtheorem*{example*}{Example}
   \newtheorem{definition}[subsubsection]{Definition}
   
   \newtheorem*{definition*}{Definition}

}
\newcommand{\CC}{{\mathbb{C}}}

\newcommand{\PP}{{\mathbb{P}}}

\renewcommand{\AA}{{\mathbb{A}}}

\newcommand{\fm}{{\boldsymbol{\mathfrak m}}}

\renewcommand{\cD}{{\mathcal D}}

\renewcommand{\cH}{{\mathcal H}}
\newcommand{\cI}{{\mathcal I}}

\newcommand{\cM}{{\mathcal M}}

\newcommand{\cO}{{\mathcal O}}
\newcommand{\cP}{{\mathcal P}}

\newcommand{\cT}{{\mathcal T}}

\def\<{\langle}
\def\>{\rangle}

\newcommand{\Spec}{\operatorname{Spec}}

\newcommand{\Proj}{\operatorname{Proj}}

\newcommand{\Hom}{{\operatorname{Hom}}}

\newcommand{\Sym}{{\operatorname{Sym}}}

\newcommand{\das}{\dashrightarrow}

\newcommand{\ocM}{\overline{{\mathcal M}}}

\newcommand{\maxord}{{\operatorname{maxord}}}

\def\:{{\colon}}
\def\.{{,\dots,}}
\def\dim{{\rm dim}}

\newcommand{\double}{\genfrac..{0pt}1
{\raise -1pt\hbox{$\scriptstyle\longrightarrow$}}{\raise 3pt\hbox
{$\scriptstyle\longrightarrow$}}}

\renewcommand{\setminus}{\smallsetminus}

\def\int{{\rm int}}

\def\tototi{\mathbin{\mathop{\otimes}\limits^{\raise-1pt\hbox
{$\scriptscriptstyle {\rm L}$}}}}

\def\indlim{\mathop{\vrule width0pt height7pt depth
4pt\smash{\lim\limits_{\raise 1pt\hbox to 14.5pt
{\rightarrowfill}}}}}
\def\projlim{\mathop{\vrule width0pt height7pt depth
4pt\smash{\lim\limits_{\raise 1pt\hbox to 14.5pt
{\leftarrowfill}}}}}

\newcommand\displaceamount{3pt}

\newcommand{\doubledown}{\ar@<\displaceamount>[d]\ar@<-\displaceamount>[d]}

\newcommand{\doubleup}{\ar@<\displaceamount>[u]\ar@<-\displaceamount>[u]}

\newcommand{\doubleright}{\ar@<\displaceamount>[r]\ar@<-\displaceamount>[r]}

\newcommand{\logord}{{\operatorname{logord}}}

\newcommand{\ord}{{\operatorname{ord}}}

\begin{document}
\title[Resolution of singularities of varieties and their families]{Resolution of singularities of complex algebraic varieties and their families}

\author{Dan Abramovich}
\address{Department of Mathematics, Box 1917, Brown University,
Providence, RI, 02912, U.S.A}
\email{abrmovic@math.brown.edu}
\thanks{New research reported here  is supported by  BSF grant 2014365, which is held jointly with J. W\l odarczyk (Purdue) and M. Temkin (Jerusalem).}

\date{\today}


\begin{abstract}
We discuss Hironaka's theorem on resolution of singularities in charactetistic 0 as well as more recent progress, both on simplifying and improving Hironaka's method of proof and on new results and directions on families of varieties, leading to joint work on toroidal orbifolds with Michael Temkin and Jaros\l aw W\l odarczyk. 
\end{abstract}
\maketitle
\setcounter{tocdepth}{1}

\tableofcontents


\section{Introduction}\label{Sec:intro}

\subsection{Varieties and singularities} An affine  complex algebraic variety $X$ is  the zero set in $\CC^n$ of a collection of polynomials  $f_i \in \CC[x_1,\ldots, x_n]$, and a general complex algebraic variety is patched together from such affine varieties much as a differentiable manifold is patched together from euclidean balls. 

But unlike differentiable manifolds, which locally are all the same (given the dimension), a complex algebraic variety can have an interesting structure locally at a point $p \in X$: the point is \emph{regular} or \emph{simple}  if the $f_i$ form the defining equations of a differentiable submanifold, and otherwise it is \emph{singular}, hiding a whole world within it. In the case of one equation, a point $p=(a_1,\ldots, a_n)$ is regular precisely when the defining equation $f_1$ has a non-vanishing derivative at $p$, and in general one needs to look at the Jacobian matrix of the defining equations, just like when studying submanifolds. The set of regular points $X^{reg} \subset X$ is always open. The variety $X$ is itself \emph{regular} if $X^{reg} = X$.\footnote{I am taking the scheme theoretic approach here: the vanishing locus of $f(x) = x^2$ is singular.}

\subsection{What is resolution of singularities?} A look at a few singularities\footnote{I cannot improve on these:\\ \url{https://imaginary.org/gallery/herwig-hauser-classic}} quickly reveals that they are quite beautiful, but complicated - really they are not \emph{simple}. How can we understand them? \emph{Resolution of singularities} provides one approach. For simplicity we restrict to \emph{irreducible} varieties, namely those which cannot be written as a union $X_1 \cup X_2$ of two closed nonempty subvarieties. A resolution of singularities of a variety $X$ is a surgery operation, a morphism $X' \to X$ which takes out the singular points and replaces them by regular points. Formally: 
\begin{definition} A \emph{resolution of singularities} of an irreducible variety $X$  is a \emph{proper morphism}  $f: X' \to X$, where $X'$  is regular and irreducible, and $f$ restricts to an isomorphism $ f^{-1} (X^{reg}) \stackrel{\sim}{\longrightarrow}  X^{reg} $.
\end{definition}

The irreducibility assumption is not serious --- for instance one can resolve each irreducible component separately.

I need to explain the terms. A \emph{morphism} $f: X' \to X$ is a mapping which locally on affine patches is defined by polynomials of the coordinates. It is an isomorphism if it is invertible as such. It is \emph{proper} if it is proper as a mapping of topological spaces in the usual Euclidean topology: the image of  a compact subset is compact. This is a way to say that we are missing no points: it would be cheating - and useless - to define $X'$ to be just $X^{reg}$: the idea is to \emph{parametrize} $X$ in a way that reveals the depths of its singularities - not to erase them!  One way $X'\to X$ can be guaranteed to be proper is if it is \emph{projective}, namely $X'$ embeds as a closed subvariety of $X \times \PP^n$ for some $n$, where $\PP^n$ stands for the complex projective space.

\subsection{Hironaka's theorem} In 1964 Hironaka published the following, see \cite[Main Theorem 1]{Hironaka}:

\begin{theorem}[Hironaka] \label{Th:Hironaka}
Let $X$ be a complex algebraic variety. Then there is a projective resolution of singularities $X' \to X$.
\end{theorem}

Hironaka's theorem is an end of an era, but also a beginning: in the half century since, people, including Hironaka, have continued to work with renewed vigor on resolution of singularities. Why is that?

I see two reasons. One reason can be seen in Grothendieck's address \cite{Grothendieck}:
\begin{quote}
Du point de vue technique, la d\'emonstration du th\'eor\`eme de Hironaka constitue
une prouesse peu commune. Le rapporteur avoue n'en avoir pas fait enti\`erement
le tour. Aboutissement d'ann\'ees d'efforts concentr\'es, elle est sans doute l'une des
d\'emonstrations les plus \guillemotleft
dures%
\guillemotright\
  et les plus monumentales qu'on connaisse en math\'ematique. 
\end{quote}
Consider, for instance, that Hironaka developed much of the theory now known as \emph{Gr\"obner bases} (at roughly the same time as Buchberger's \cite{Buchberger-thesis}) for the  purpose of resolution of singularities! 

There has been a monumental effort indeed to simplify Hironaka's proof, and to break it down to more basic  elements, so that the techniques involved come naturally and the ideas flow without undue effort. I think this has been a resounding success and Grothendieck himself would have approved of the current versions of the proof - he certainly would no longer have trouble going through it. In my exposition I attempt to broadly describe the results of this effort.

A few points in this effort  are marked by the following: 
\begin{itemize}
\item The theory of maximal contact, \cite{Giraud,AHV}.
\item Constructive resolution using an invariant, \cite{Villamayor, Bierstone-Milman, Encinas-Villamayor-good}.
\item The optimal version of canonical resolution \cite{Bierstone-Milman}.
\item Simplification using order of ideal \cite{Encinas-Villamayor}.
\item Functoriality as a proof technique and guiding principle  \cite{Wlodarczyk, Bierstone-Milman-funct}.
\item Dissemination to the masses \cite{Cutkosky-book, Hauser-stories, Kollar}.
\end{itemize}

The other reason is generalizations and refinements of the resolution theorem. First and foremost, algebraic geometers want to resolve singularities in positive and mixed characteristics, as the implications would be immense.  In addition, one is interested in simplifying families of varieties, simplifying algebraic differential equations, making the resolution process as effective and as canonical as possible, and preserving structure one is provided with at the outset. Below I will discuss  \emph{resolution in families} while preserving \emph{toroidal} structures, focussing on joint work with Michael Temkin of Jerusalem and Jaros\l aw W\l odarczyk of Purdue.

\section{Hironaka's method: from resolution to order reduction}

The purpose of this section is to indicate how Hironaka's resolution of singularities can be reduced to an algebraic problem, namely \emph{order reduction} of an ideal.

\subsection{Blowing up} The key tool for Hironaka's resolution of singularities is an operation called \emph{blowing up} of a regular subvariety $Z$ of a regular variety $Y$, see \cite[Definition p. 163]{Hartshorne}

\subsubsection{Blowing up a point} A good idea can be gleaned from the special case where $Y = \AA^n$, affine $n$-space, and $Z$ is the origin, as explained in \cite[Example 7.12.1]{Hartshorne} and depicted on the cover of \cite{Shafarevich1}. Think about $\PP^{n-1}$ as the set of lines in $\AA^n$ through the origin. The blowing up $Y' \to Y$ is then given as the \emph{incidence variety} 
$$Y' = \{ (x,\ell) \in \AA^n \times \PP^{n-1} \mid x \in \ell\}, $$ with its natural projection to $\AA^n$.

This can be described in equations as follows:
$$Y' = \{ ((x_1,\ldots, x_n),(Y_1:\,\ldots\,: Y_n)) \in \AA^n \times \PP^{n-1} \mid x_i Y_j = x_j Y_i\ \forall  i,j\}.$$

Since $\PP^n$ is covered by affine charts, this can further be simplified. For instance on the chart where $Y_n \neq 0$ with coordinates $y_1,\ldots, y_{n-1}$, this translates to
$$\{ ((x_1,\ldots, x_n),(y_1,\ldots, y_{n-1})) \in \AA^n \times \AA^{n-1} \mid x_j = x_n y_j, \ 1\leq j\leq n-1 \}.$$ In other words, the coordinates $x_1,\ldots,x_{n-1}$ are redundant and this  is just affine space with coordinates $y_1,\ldots, y_{n-1},x_n$. In terms of coordinates on $Y$ we have $y_j = x_j/x_n$.

The fibers of $Y' \to Y$ are easy to describe: away from the origin $x_1 = \dots=x_n = 0$ the map is invertible, as the line $\ell$ is uniquely determined by $(x_1,\dots, x_n)$. Over the origin all possible lines occur, so the fiber is $\PP^{n-1}$, naturally identified as the space of lines through the origin.

\subsubsection{Blowing up a regular subvariety} In general the process is similar: given regular subvariety $Z$ of $Y$, then $f:Y'\to Y$ replaces each point $z\in Z$ by the projective space of normal directions to $Z$ at $z$. If $Z$  is locally defined by equations $x_1=\ldots=x_k=0$ and if  $x_{k+1}, \ldots, x_n$ form coordinates along $Z$, then $Y'$ has local patches, one corresponding to each $x_i$, with coordinates 
$$\frac{x_1}{x_i}, \ldots, \frac{x_{i-1}}{x_i}, \ x_i,\  \frac{x_{i+1}}{x_i},\ldots, \frac{x_k}{x_i},\  x_{k+1}, \ldots, x_n.$$

Thus the blowing up $Y'\to Y$ of a regular subvariety $Z$ of a regular variety $Y$ always results in a regular variety.

We often say that $Z$ is the \emph{center} of the blowing up $Y'\to Y$, or that the blowing up $Y'\to Y$ is \emph{centered} at $Z$.

\subsubsection{The $\Proj$ construction}
Grothendieck gave a more conceptual construction, which applies to an arbitrary subscheme $Z$ defined by an ideal sheaf $\cI$ in an arbitrary scheme $Y$:

$$Y' = \Proj_Y \bigoplus_{k=0}^\infty \cI^k / \cI^{k+1}.$$

The map $f^{-1} (Y\setminus Z) \to (Y\setminus Z)$ is always an isomorphism, so we identify $Y\setminus Z$ with its preimage.  

What the reader may want to take from this is that the blowings up we introduced in particularly nice cases are part of a flexible array of transformations.

The complement $E$ of $Y\setminus Z$ in $Y'$ is called the \emph{exceptional locus}. It is  \emph{a Cartier divisor}, a subvariety of  codimension 1 locally defined by one equation. If $Z$ is nowhere dense in $Y$, then $Y' \to Y$ is birational. If moreover $Y$ and $Z$ are regular, then $E$ is regular. 

\subsubsection{The strict transform}

Blowing up serves in the resolution of singularities of a subvariety $X \subset Y$ through the \emph{strict transform} $X' \subset Y'$: this is the closure of $X\setminus Z$ in $Y'$. Grothendieck showed that $X'$ is the same as the blowing up of $X\cap Z$ in $X$, using the $\Proj$ construction above.

From the point of view of resolution of singularities, the challenge is to make $X'$ less singular than $X$ by an appropriate choice of $Z$. 

Consider for instance the cuspidal plane curve $X$ given by $y^2 - x^3=0$ in the affine plane $Y=\AA^2$ with coordinates $x,y$. Blowing up the origin and focusing on the chart with coordinates $x,z=y/x$, we obtain the equation $z^2 x^2 - x^3=0$, which we rewrite as $x^2 (z^2-x)=0$. The locus $x=0$ describes the exceptional line, and $X'$  is given by $z^2 - x=0$, a regular curve.

\subsection{Embedded resolution}

Theorem \ref{Th:Hironaka} is proven by way of the following theorem:

\begin{theorem}[Embedded resolution]\label{Th:embedded} Suppose $X \subset Y$ is a closed subvariety of a regular variety $Y$. There is a sequence of blowings up $Y_n \to Y_{n-1} \to \cdots \to Y_0 = Y$, with regular centers $Z_i \subset Y_i$ and strict transforms $X_i \subset Y_i$, such that $Z_i$ does not contain any irreducible component of $X_i$  and such that $X_n$ is regular.
\end{theorem}

In other words, the final strict transform of $X$ in a suitably chosen  sequence of blowings up of $Y$ is regular.

In the example of the cuspidal curve $X \subset Y=\AA^2$, the single blowing up $Y_1 \to Y$ centered at the origin $Z = \{(0,0)\}$ provides an embedded resolution $X' \to X$ of $X$.

If $X$ is embedded inside a regular variety $Y$,  then Theorem \ref{Th:embedded} immediately gives a resolution of singularities  $X_n \to X$. What if $X$ is not embedded? There are a number of viable approaches, but the best is to strengthen Theorem \ref{Th:embedded}: one makes the blowing up procedure independent of re-embedding $X$ and compatible with local patching. This is what is done in practice. I describe the underlying principles in Sections \ref{Sec:functoriality} and \ref{Sec:reembedding} below. The upshot is that ``good embedded resolution implies resolution''. 

From here on we pursue a good embedded resolution.

\subsection{Normal crossings} To go further one needs to describe a desirable property of the exceptional divisor $E_i$ and its interaction with the center $Z_i$. 

\begin{definition} We say that a closed subset $E \subset Y$ of a regular variety $Y$ is a \emph{simple normal crossings divisor} if in its decomposition $E = \cup E_j$ into irreducible components, each component $E_j$ is regular, and these components intersect transversally: locally at a point $p\in E$ there are local parameters $x_1,\ldots,x_m$ such that $E$ is the zero locus of a reduced monomial $x_1\cdots x_k$.

We further say that $E$ and a regular subvariety $Z$ \emph{have normal crossings} if such coordinates can be chosen so that $Z=V(x_{j_1}, \ldots, x_{j_l})$ is the zero set of a subset of these coordinates.
\end{definition} 

When the set of coordinates $x_{j_1}, \ldots, x_{j_l}$ is disjoint from $x_1,\dots,x_k$ the strata of $E$ meet $Z$ transversely, but the definition above allows quite a bit more flexibility.

This definition works well with blowing up: If $E$ is a simple normal crossings divisor,  $E$ and $Z$ have normal crossings,  $f:Y'\to Y$ is the blowing up of the regular center $Z$ with exceptional divisor $E_Z$, and $E' = f^{-1}E \cup E_Z$ then $E'$ is a simple normal crossings divisor.

\subsection{Principalization} Embedded resolution is proven by way of the following algebraic result:

\begin{theorem}[Principalization] Let $Y$ be a regular variety and $\cI$ an ideal sheaf. There is a sequence of blowings up $Y_n \to Y_{n-1} \to \cdots \to Y_0 = Y$, regular subvarieties $Z_i\subset Y_i, i=0,\ldots,n-1$ and simple normal crossings divisors $E_i \subset Y_i, i=1,\ldots,n $ such that
\begin{itemize}
\item $f_i:Y_{i+1} \to Y_i$ is the blowing up of $Z_i$ for $ i=0,\ldots,n-1$,
\item $E_i$ and $Z_i$ have normal crossings for $i=1,\ldots,n-1$,
\item $\cI\cO_{Y_i}$ vanishes on $Z_i$ for $ i=0,\ldots,n-1$, 
\item $E_{i+1}$ is the union of $f_i^{-1} E_i$ with the exceptional locus of $f_i$ for $ i=0,\ldots,n-1$
\end{itemize}
 and such that the resulting ideal sheaf $\cI_n = \cI \cO_{Y_n}$ is an invertible ideal with zero set $V(\cI_n)$ supported in $E_n$.
\end{theorem}

In local coordinates $x_1,\ldots,x_m$ on $Y_n$ as above, this means that $\cI_n = (x_1^{a_1} \cdots x_m^{a_k})$ is locally principal and monomial, hence the name ``principalization''. The condition that $Z_i$ have  normal crossings with $E_i$ guarantees that $E_{i+1}$ is a simple normal crossings divisor.

\subsubsection{Principalization implies embedded resolution} Quoting Koll\'ar \cite[p. 137]{Kollar}, principalization implies embedded resolution seemingly ``by accident'': suppose for simplicity that $X$ is irreducible, and let the ideal of $X\subset Y$ be $\cI$. Since $\cI_n$ is the ideal of a divisor supported in the exceptional locus, at some point in the sequence the center $Z_i$ must contain the strict transform $X_i$ of $X$. Since $\cI$ vanishes on $Z_i$, it follows that $Z_i$ coincides with $X_i$ at least near $X_i$. In particular $X_i$ is regular! 

\subsubsection{Are we working too hard?} Principalization  seems to require ``too much'' for resolution: why should we care about exceptional divisors which lie outside $X$? Are we trying too hard?

In the example of the cuspidal curve $X \subset Y = \AA^2$ above, the single blowing up $Y_1 \to Y$ at the origin does \emph{not} suffice for principalization: the resulting equation $x^2(z^2-x)=0$ with exceptional $\{x=0\}$ is \emph{not} monomial. One needs no less that \emph{three} more blowings up! I'll describe just one key affine patch of each:
\begin{itemize}
\item Blowing up $x=z=0$ one gets, in one affine patch where $x=zw$, the  equation $z^3 w^2 (z-w)=0$, strict transform $X_2= \{z=w\}$ and exceptional $\{wz=0\}$. 
\item Blowing up $z=w=0$ one gets, in one affine patch where $z=wv$, the equation  $v^3 w^6 (v-1)=0$. The exceptional in this patch is $\{wv=0\}$, with one component (the old $\{z=0\}$) appearing only in the other patch. The strict transform is $X_3 = \{v=1\}$.
\item 
In the open set $\{v\neq 0\}$ we blow up $\{v=1\}$. This actually does nothing, except turning the function $u=v-1$ into a monomial along $v=1$, so the equation $v^3 w^6 (v-1)=0$ at these points can be written as $(u+1)^3w^6u=0$, which in this patch  is equivalent to $w^6u=0$, a monomial in the exceptional parameters $u,w$.
\end{itemize}

The fact that we could blow up $\{v=1\}$ means that $X_3$ is regular, giving rather late evidence that we obtained resolution of singularities for $X$. These ``redundant'' steps add to the sense that this method works  ``by accident''. It turns out that principalization itself is quite useful in the study of singularities. Also the fact that it provides the prize of resolution is seen as sufficient justification.  The discussion in Section \ref{Sec:ATW} will put it in the natural general framework of toroidal structures.

Accident or not, we will continue to pursue principalization.

\subsection{Order reduction}\label{Sec:order-reduction} Finally, principalization of an ideal is proven by way of \emph{order reduction}.

The \emph{order} $\ord_p (\cI)$ of an ideal $\cI$ at a point $p$ of  a regular variety $Y$ is the maximum integer $d$ such that $\fm_p^d \supseteq \cI$; here $\fm_p$ is the maximal ideal of $p$. It tells us ``how many times  every element of $\cI_p$ vanishes at $p$.'' 

In particular we have $\ord_p(\cI) \geq 1$ precisely if $\cI$ vanishes at $p$. 

We write $\maxord(\cI) = \max\{\ord_p(\cI) | p\in X\}$. For instance we have $\maxord(\cI) = 0$ if and only if $\cI$ is the unit ideal, which vanishes nowhere. Another exceptional case is $\maxord(\cI)=\infty$ which happens if $\cI$ vanishes on a whole component of $Y$. We'll ignore that case for now.

Given an integer $a$, we write $V(\cI,a)$ for the locus of points $p$ where $\ord_p(\cI) \geq a$. A regular closed subvariety $Z\subset Y$ is said to be \emph{$(\cI,a)$-admissible} if and only if $Z \subset V(\cI,a)$, in other words, the order of $\cI$ at every point of $Z$ is at least $a$. Admissibility is related to blowings up: if $\maxord(\cI) =a$, and if $Y' \to Y$ is the blowing up of an $(\cI,a)$-admissible $Z\subset Y$, with exceptional divisor $E$ having ideal $\cI_E$, then $\cI \cO_{Y'} = (\cI_E)^a\, \cI'$, with $\maxord(\cI') \leq a$.

\emph{Order reduction} is the following statement:

\begin{theorem}[Order reduction]\label{Th:order-reduction} Let $Y$ be a regular variety, $E_0\subset Y$ a simple normal crossings divisor, and $\cI$ an ideal sheaf, with $$\maxord(\cI) = a.$$ There is a sequence of $(\cI,a)$-admissible blowings up $Y_n \to Y_{n-1} \to \cdots \to Y_0 = Y$, with regular centers $Z_i \subset Y_i$ having normal crossings with $E_i$ such that  $ \cI \cO_{X_n} = \cI_n \cI_n'$
with $\cI_n$ an invertible ideal supported on $E_n$ and such that  $$\maxord(\cI_n')<a.$$
\end{theorem}

Order reduction implies principalization simply by induction on the maximal order $\maxord(\cI) =a$: once $\maxord(\cI_n')=0$ we have $ \cI \cO_{X_n} = \cI_n$ so only the exceptional part remains, which is supported on a simple normal crossings divisor by induction.

Hironaka himself used the \emph{Hilbert--Samuel function}, an invariant much more refined than the order. It is a surprising phenomenon that resolution becomes easier to explain when one uses just the order, thus less information, see \cite{Encinas-Villamayor}.

\emph{It remains to prove order reduction.}

\section{Hironaka's method: order reduction}

\subsection{Differential operators} Nothing so far was particularly sensitive to the fact that we were working over $\CC$, or even a field of characteristic 0. That starts changing now.

Since $Y$ is regular, it has a \emph{tangent bundle} $T_Y$. Local sections $\partial$ of $T_Y$ are first order differential operators $\partial: \cO_Y \to \cO_Y$. As usual we denote the sheaf of sections of the tangent bundle with the same symbol $T_Y$, hoping the confusion can be overcome. 

\subsubsection{The characteristic 0 case} In characteristic 0, the sheaf of rings generated over $\cO_Y$ by the operators in $T_Y$ is the \emph{sheaf of differential operators} $\cD_Y$. As a sheaf of $\cO_Y$ modules it looks locally like the symmetric algebra $\Sym^\bullet(T_Y) = \oplus_{n\geq 0} \Sym^n(T_Y)$, but its ring structure is very different, as  $\cD_Y$ is non-commutative. Still for any integer $a$ there is a subsheaf $\cD^{\leq a}_Y \subset \cD_Y$ of differential operators of order $\leq a$, those sections which can be written in terms of monomials of order at most $a$ in sections of $T_Y$. As a special case, one always has a splitting $\cD^{\leq 1}_Y = \cO_Y \oplus T_Y$, the projection $\cD^{\leq 1}_Y \to \cO_Y$ given by applying $\nabla \mapsto \nabla(1)$. 

\subsubsection{The general case}  Things are quite different in characteristic $p>0$: one can use the same definition, but in some sense it is deficient, because these differential operators do not detect $p$th powers. There is a natural and sophisticated replacement, which coincides with $\cD_Y$ in characteristic $0$, and defined as follows:

On $Y\times Y$ consider the diagonal $\Delta \subset Y\times Y$. It is a closed subvariety, and one can consider its ideal $\cI_\Delta$. The \emph{sheaf of principal parts of order $a$ of $Y$} is defined as $\cP\cP^a_Y = \cO_{Y\times Y} / \cI_\Delta^{a+1}$ - it is a sheaf of $\cO_Y$-modules via either projection; its fiber at $p\in Y$ describes functions on $Y$ up to order $a$ at $p$. Its dual sheaf is $\cD^{\leq a} := (\cP\cP^a_Y)^\vee$, which \emph{in characteristic 0} admits the concrete description given earlier. The natural projection $ \cP\cP^a_Y \to \cP\cP^{a-1}_Y$ gives rise to an inclusion $\cD_Y^{\leq a-1} \subset \cD_Y^{\leq a}$, and one defines in general $\cD_Y = \cup_a \cD^{\leq a}_Y$.

This is nice enough, but the fact that in positive characteristics sections of $\cD_Y$ are not written as polynomials in sections of $T_Y$ is the source of much trouble.

\subsection{Derivatives and order}

Let $\cI$ be an ideal sheaf on $Y$ and $y\in Y$ a point. Write $\cD_Y^{\leq a}\cI$ for the ideal generated by elements $\nabla (f)$, where $\nabla$ an operator in  $\cD_Y^{\leq a}$ and $f$ a section of $\cI$. We have the following characterization:

$$ \ord_y(\cI) = \min\{a: (\cD^{\leq a} \cI)_p = \cO_{Y,p} \}.$$

In other words, the order of $\cI$ at $y$ is the minimum order of a differential operator $\nabla$ such that for some $f \in \cI_y$  the element $\nabla(f)$ does not vanish at $y$. 

We define $\cT(\cI,a):= \cD^{\leq a-1} \cI$. In these terms, the set $V(\cI,a)$ can be promoted to a \emph{scheme}, the zero locus of an ideal: $V(\cI,a) := V(\cT(\cI,a))$.\footnote{This is the \emph{right} scheme structure, as it satisfies an appropriate universal property.} 

This is not too surprising in characteristic $0$ since we all learned calculus, but it may seem strange in characteristic $p>0$. For instance, the order of $(x^p)$ is $p$, since there is always an operator $\nabla$ of order $p$ such that $\nabla (x^p) = 1$. In characteristic $0$ we can write $$\nabla = \frac1{p!} \left(\frac{\partial}{\partial x}\right)^p,$$ but in characteristic $p$ we have no such expression! 

\subsection{Induction and maximal contact hypersurfaces} We return to working over $\CC$, in particular in characteristic 0, so we can use the letter $p$ for a point of $Y$.

Remember that we want to prove order reduction of an ideal $\cI$ of maximal order $a$. Hironaka's next idea was to use induction on dimension by restricting attention to a hypersurface $H$, in such a way that a suitable order reduction on $H$ results, \emph{by blowing up the same centers,} in order reduction of $\cI$ on $Y$. 

I am not being historically correct here, since Hironaka used invariants much more refined than order. I depart from history further, and use Giraud's concept of \emph{maximal contact hypersurfaces}, adapted to orders rather than other invariants.  

\begin{definition} Let $\cI$ be an ideal of maximal order $a$. \emph{A maximal contact hypersurface for $(\cI,a)$ at $p$} is a hypersurface $H$ \emph{regular at $p$}, such that, in some neighborhood $Y^0$ of $p$ we have $H \supseteq V(\cI,a) = V(\cT(\cI,a))$, namely $H$ contains  the \emph{scheme} of points where $\cI$ has order $a$.
\end{definition}

\subsection{Derivatives and existence in characteristic 0} It is not too difficult to show that \emph{in characteristic 0}, a maximal contact hypersurface for $(\cI,a)$ at $p$ exists. Since $\cI$ has maximal order $a$, we have $\cD^{\leq a}(\cI) = (1)$. Consider the ideal $\cT(\cI,a)=\cD^{\leq a-1} (\cI)$. Since we are in characteristic 0, it must contain an antiderivative of $1$, so $\ord_p(\cT(\cI,a)) \leq 1$. Any local section $x$ of $\cT(\cI,a)$ with order $\leq 1$ gives a maximal contact hypersurface $\{x=0\}$ at $p$.

Here is an example: suppose $\cI= (f)$ where $$f(x,y) = y^a + g_1(x)y^{a-1} + \dots + g_{a-1}(x)y+ g_a(x).$$ Then $\ord_{(0,0)}(f) = a$ exactly when $\ord_0(g_i) \geq i$ for all $i$. In characteristic $0$ we may replace $y$ by $y+g_1(x)/a$, so we may assume $g_1(x) = 0$. In this case $\partial^{a-1}f / \partial y^{a-1} = a!\cdot y$, so $\{y=0\}$ is a maximal contact hypersurface at $(0,0)$.

The definition I gave is pointwise. It is easy to see that if $H$ is a maximal contact hypersurface for $(\cI,a)$ at $p$ then the same holds at any nearby $p'$, so the concept is local. Unfortunately it is also not hard to cook up examples where there is no \emph{global} maximal contact hypersurface which works everywhere. We will tackle this problem in section \ref{Sec:functoriality} below.

\subsubsection{Positive characteristics} Alas, there are fairly simple examples in charactersitic $p>0$ where maximal contact hypersurfaces do not exist, \cite{Narasimhan, Hauser-positive}. The whole discussion from here on simply does not work in characteristic $>0$.

\subsection{What should we resolve on $H$?}
For induction to work we need to decide what exactly we want to do on the hypersurface $H$. The example $\cI= (f)$ above is instructive: just restricting $\cI$ to $H$ does not work! 

Assume $$f(x,y) = y^a + g_2(x)y^{a-2} + \dots + g_{a-1}(x)y+ g_a(x),$$ with $\ord_0(g_i) \geq i$ for all $i$. The restriction of $\cI$ to the hypersurface $H = \{y=0\}$ is the ideal $(g_a(x))$. Now clearly this ideal does not retain enough information from the original ideal. This is manifest with the notion of admissibility introduced in Section \ref{Sec:order-reduction}, as   $(g_a(x),a)$-admissible centers in $H$ will not always give $(\cI,a)$-admissible centers in $Y$. For instance it might happen that $g_a=0$, so every center on $H$ is $(g_a(x),a)$-admissible, but $g_{a-1}\neq 0$, so not every center on $H$ is $(\cI,a)$-admissible on $Y$!

The collection of elements $g_2(x),\ldots, g_a(x)$ surely hold all the necessary information. However each comes with its own requirements: in order to reduce the order of $\cI$ below $a$, we need to reduce the order of at least one of $g_i(x)$ below $i$. 

We now generalize this discussion to arbitrary ideals.

\subsubsection{Coefficient ideals and the induction scheme}
In order to generalize this, we need to identify an analogue of these ``elements'' $g_i(x)$, and derivatives come to the rescue again. Let $\cI$ be an ideal of maximal order $a$ on a variety $Y$ and $H$ a maximal contact hypersurface. Then for any $i< a$ the ideals $\cD_Y^{\leq i} \cI$ have maximal order precisely $a-i$. It follows that the restricted ideals $(\cD_Y^{\leq i} \cI)_{|H}$ have maximal order $\geq a-i$. These restrictions are the analogues of $g_{a-i}(x)$. 

The following is at the technical core of the proof. I am aware of several proofs of this proposition, but they all seem to go a bit beyond the level of discussion I wish to maintain here. For reference, see  \cite[Section 3.9]{Kollar}.

\begin{proposition} Any sequence of $(\cI,a)$-admissible blowings up has centers lying in $H$ and its successive strict transforms. The resulting sequence of blowings up on $H$ is $((\cD_Y^{\leq i} \cI)_{|H}, a-i)$-admissible for every $i<a$. 

Conversely, every sequence of blowings up on $H$ which is $((\cD_Y^{\leq i} \cI)_{|H}, a-i)$-admissible for every $i<a$ gives rise, by blowing up the same centers on $Y$, to a sequence of $(\cI,a)$-admissible blowings up.
\end{proposition}

Such an admissible sequence on $H$ may be called an \emph{order reduction} for the collection $((\cD_Y^{\leq i} \cI)_{|H}, a-i)$ if it forms an order reduction for at least one of these pairs. It is a formal consequence of the proposition that an order reduction for $(\cI,a)$ is the same as an order reduction for the collection $((\cD_Y^{\leq i} \cI)_{|H}, a-i), i=0,\ldots, a.$

This may appear as troublesome: we wanted to prove order reduction for one ideal, and the induction requires us to prove order reduction for a collection of ideals. But there is a simple trick that allows one to replace this collection of ideals by a single ideal, in such a way that the notions of order reduction coincide: 

\begin{definition} For an ideal $\cI$ of maximal order $a$ define its \emph{coefficient ideal} to be the ideal sum
$C(\cI,a) :=\sum (\cD_Y^{\leq i} \cI)^{a!/(a-i)}.$
\end{definition}
\begin{proposition}[{\cite[\S 3.4]{Wlodarczyk}}]\label{Prop:lifting} Order reduction for $(C(\cI,a)_{|H},a!)$ is the same as order reduction for the collection $((\cD_Y^{\leq i} \cI)_{|H}, a-i), i=0,\ldots, a.$
\end{proposition}

We obtain:
\begin{corollary}[{\cite[Corollary 3.85]{Kollar}}]\label{Cor:lifting} A sequence of blowings up is an order reduction for $(\cI,a)$ if and only if it is an order reduction for $(C(\cI,a)_{|H}, a!)$.
\end{corollary}

\subsection{A trouble of exceptional loci} I have been deliberately ignoring a subtle point. Theorem \ref{Th:order-reduction} about order reduction takes the additional  datum of a divisor $E_0 \subset Y$. This is important for principalization, since once we reduce the order of $\cI$ from $a$ to $a-1$ with exceptional divisor say $E_0$, we want any further centers of blowing up used in further order reduction of $\cI$ to have normal crossings with $E_0$.

For instance, in the example of a cuspidal curve above, the ideal $\cI= x^2 (z^2-x)$ is of the form $\cI_E^2 \cI'$. The unique maximal contact hypersurface for $(\cI',1)$ is precisely $X'$, the vanishing locus of $\cI'$, but since it is tangent to $E$ it  does not have normal crossings with $E$. 

The standard way to treat this is via a trick: one separates the relevant part of the ideal $\cI'$ from the monomial part $\cI_E$ by applying a suitable principalization for an ideal of the form $\cI_E^\alpha + {\cI'}^\beta$ describing the intersection of their loci. This is somewhat subtle and a bit disappointing. One feels that monomial ideals should only serve for good, as they are the goal. 

I'll totally ignore this issue here, referring to \cite[Section 3.13]{Kollar}. I have an excuse: in the procedure described below in my work with Temkin and W\l odarczyk, this is not an issue at all, as monomial ideals become our best friends. 

\subsection{The problem of gluing}

I postponed two important issues. Resolution of singularities requires \emph{good} embedded resolution, \emph{good} principalization, \emph{good} order reduction:  the process must be compatible with  patching of open sets and independent of the embedding. A related issue is the fact that maximal contact hypersurfaces  are not global, so patching open sets where maximal contact hypersurfaces do not overlap is required!

The classical approach has several ideas involved and has several levels of complexity. First, one devises a more elaborate \emph{resolution invariant}, which records behavior of a given ideal on a sequence of nested maximal contact hypersurfaces. Second, one devises a class of transformations, called \emph{test transformations},  which include admissible blowings up, restrictions to open sets,  but also other operations, such as projections from a product. 

One shows that ideals with the same invariant admit the same sequences of test transformations. It is a more subtle fact that the opposite is true - the invariant can be read off the test transformation. Once the dust settles it becomes clear that the order reduction one produces is independent of choices and is local, hence it can be patched along open sets. Also, the issue of the choice of embedding for resolution of singularities becomes local, hence it reduces to a simple principle I call \emph{the re-embedding principle} in Section \ref{Sec:reembedding} below.

Some subset of this approach, in particular the fact that invariants can be read from the class of test transformations, is known as \emph{Hironaka's trick}.

I have to admit that I never quite understood this approach until I read W\l odarczyk's paper \cite{Wlodarczyk}, which uses a completely different approach. After that transformative event I was able to read Bierston and Milman's \cite{Bierstone-Milman-funct}, and suddenly the classical approach was illuminated. I therefore prefer to present W\l odarczyk's approach here, as perhaps  others will experience the same transformation and subsequent illumination.

\subsection{W\l odarczyk's functoriality principle}\label{Sec:functoriality}

Already Hironaka was interested in functorial properties of resolution of singularities. For instance, if $X$ is a singular variety with a group $G$ acting, ideally one would want the resolution to be $G$ equivariant. Also if $X^0 \subset X$ is open, the resolution of $X^0$ should ideally be the restriction of that of $X$. This is stated explicitly in \cite[\S13]{Bierstone-Milman}.

W\l odarczyk's great idea in \cite{Wlodarczyk} was that \begin{quote} \em functoriality is a powerful tool in the very proof of resolution of singularities. \end{quote}
Moreover, \begin{quote} \em functoriality leads one to discover an order reduction algorithm. \end{quote}
W\l odarcyk requires one to take this very seriously. Indeed, for his principle to succeed one needs to use hidden symmetries, which are revealed only after $\cI$ is tuned appropriately.

\subsubsection{Smooth pullbacks} Let us first define the terms. Let $Y_n \to \dots \to Y_0=Y$ be an order reduction of $(\cI,a)$ compatible with simple normal crossings divisor $E$. Let $Y'\to Y$ be a \emph{smooth} morphism, what geometers call a \emph{submersion}, such as an open embedding or a product with a regular variety. One can write $\cI' = \cI\cO_Y$, $E'=E\times_Y Y'$ and $Y_i' = Y_i\times_Y Y'$, and then automatically  $Y'_n \to \dots \to Y'_0=Y'$ is an order reduction for $(\cI', a)$, compatible with $E'$, the \emph{smooth pullback} order reduction. Some of the resulting steps might become trivial, in which case we drop them from the order reduction sequence.

\begin{definition} A \emph{functorial order reduction} is a rule assigning to an ideal $\cI$ on a regular variety $Y$ with simple normal crossings divisor $E$  and integer $a$ such that $\maxord(\cI) \leq a$, an order reduction  $Y_n\to \cdots \to Y$, in such a way that for any smooth morphism $Y' \to Y$, the corresponding order reduction $Y_n' \to \cdots \to Y'$ is the smooth pullback order reduction of $Y_n\to \cdots \to Y$.
\end{definition}

\emph{To prove order reduction it suffices to produce \emph{functorial} order reduction on open patches, because then they automatically glue together.} 

What about maximal contact hypersurfaces? Let us say we have produced functorial order reduction in dimension $\dim (Y) - 1$ and we wish to prove it for $Y$. We can choose a local maximal contact hypersurface $H \subset Y$ and reduce the order of the coefficient ideal $C(\cI,a)_{|H}$. By Corollary \ref{Cor:lifting} this results in \emph{local} order reduction for $(\cI,a)$, but \emph{a priori} this depends on the choice of $H$. We claim that in fact there is no such dependence, and that the resulting order reduction is functorial on $Y$. For this we use $(\cI,a)$-special automorphisms.

\subsubsection{Special automorphisms} Let $Y_m \to \dots \to Y_0$ be an $(\cI,a)$-admissible sequence with centers $Z_i \subset Y_i$. Recall that this gives in  particular  a sequence of ideals $\cI_i \subset \cO_{Y_i}$ such that $Z_i \subset V(\cI_i,a)$ and $\cI_i \cO_{y_{i+1}} = \cI_{E_{i+1}}^a \cI_{i+1}$. An automorphism $\phi$ of $Y$ is \emph{special} if it \emph{fixes} every $(\cI,a)$-admissible sequence. This means that $\phi$ fixes  $V(\cI,a)$, in particular it fixes $Z_0$, hence it lifts to an automorphism $\phi_1$ of $Y_1$, \emph{which fixes $V(\cI_1,a)$}, and inductively we obtain automorphisms $\phi_i$ of $Y_i$ fixing $V(\cI_i,a)$. 

This is a very strong assumption on an automorphism, but W\l odarczyk proved the following powerful result:

\begin{proposition}[{\cite{Wlodarczyk}}]\label{Prop:loca-aut} Let $H_1,H_2\subset Y$ be two local maximal contact hypersurfaces at $p\in V(\cI,a)$. Then, after replacing $Y$ by an \'etale neighborhood of $p$, there is a special automorphism $\phi$ of $Y$ fixing $p$ and sending $H_1$ to $H_2$.
\end{proposition}

In particular, the functorial order reductions for $C(\cI,a)_{|H_i}$ induce the same  order reduction for $\cI$, which is automatically functorial!

I deliberately did not require the automorphism $\phi$ to send $\cI_i$ to itself, which would make the statement easier to grasp. Indeed,  
the following example shows that in general \emph{it is impossible for $\phi$ to send $\cI_i$ to itself,} and suggests that Proposition \ref{Prop:loca-aut} is quite surprising and should require an ingenious idea. 

Consider the ideal $(xy)$ in the affine plane $Y$, with maximal order $2$ attained at $V(\cI,2) = (0,0)$, the origin. We have $\cD^{\leq 2-1}\cI = \cD^{\leq 1} \cI = (x,y)$, and so the lines $H_1 = \{x=0\}$ and $H_2 = \{x+y=0\}$ are both maximal contact hypersurfaces. Clearly any automorphism of $Y$ sending $H_1$ to $H_2$ must change $\cI$, since $\cI_{|H_1} =0$ and $\cI_{|H_2} \neq 0$. 

In this particular case the coefficient ideal is $(x^2, xy, y^2)$, and the automorphism $(x,y) \mapsto (x+y, y)$ does send $H_1$ to $H_2$ fixing this ideal, so whatever procedure we apply using $H_1$  - in this case necessarily blowing up the origin - coincides with the process we apply using $H_2$.

\subsubsection{Homogenization} The general case is slightly more subtle than the example: in general there is no automorphism carrying $H_1$ to $H_2$ fixing the coefficient ideal either. Searching for a natural replacement which is fixed under a special automorphism, W\l odarcyk discovered the \emph{homogenization} $\cH(\cI,a)$ described below.\footnote{A Different variant is used in \cite{Kollar}; the treatment in  \cite{BGV} in terms of differential Rees algebras provides a natural structure subsuming homogenization and coefficient ideals.}

\begin{definition} Recall the notation  $\cT(\cI,a) = \cD^{\leq a-1}_Y \cI$. The \emph{homogenization} of $(\cI,a)$ is the ideal
$$\cH(\cI, a) := \sum_{i=0}^{a}\cD^i(\cI)\, \cT(\cI,a)^i.$$ 
\end{definition}

W\l odarczyk's idea is that $\cI$ lacks symmetries because it is not sufficiently tuned. In contrast, the ideal $\cH(\cI,a)$ is tuned to reveal the hidden symmetry $\phi$.\footnote{In \cite{Bierstone-Milman-funct}, Bierstone and Milman replace $(\cI,a)$ by its equivalence class with respect to test transformation. With the ``blurred vision'' of equivalence classes of ideals, a hidden symmetry is again revealed. This is related to Hironaka's approach using the concept of \emph{infinitely near points}.}

The ideal $\cH(\cI,a)$ is designed to contain all terms of Taylor expansions of elements of $\cI$ in terms of any variable $h$ in $\cT(\cI,a)$. If $H_1 = \{x=0\}$, $H_2 = \{x+h = 0\}$ and $x=x_1, x_2,\ldots,x_m$ are local parameters of $Y$ and $p$, chosen so that  $x_1+h, x_2,\ldots,x_m$ also form local parameters, then the transformation $\phi(x_1,x_2,\ldots,x_m)=(x_1+h,x_2,\ldots,x_m)$ is a local automorphism of  $Y$ formally sending $f(x_1,x_2,\ldots,x_m)$ to $$\sum \frac{\partial^if}{\partial x_1^i} h^i.$$ Note that $\frac{\partial^if}{\partial x_1^i} h^i \in \cD^i(\cI)\, \cT(\cI,a)^i$. Thus on formal completions this sends an element of $\cH(\cI,a)$ to an element of $\cH(\cI,a)$, and a bit of reflection shows that $\phi$ is a special automorphism with respect to $\cH(\cI,a)$. A standard argument allows to pass from completion to \'etale neighborhoods, hence $\phi$ defines a special automorphism with respect to $\cH(\cI,a)$ on a suitable \'etale neighborhood. 

A simple computation shows:

\begin{proposition} Order reduction for $(\cI,a)$ is equivalent to order reduction for $\cH(\cI,a)$.
\end{proposition}

Thus $\phi$ is a special automorphism with respect to $(\cI,a)$ as well!

\subsection{A sketch of the algorithm} Let us summarize how one \emph{functorially} reduces the order of a nonzero ideal $\cI$ of maximal order $a>0$ on a regular variety $Y$. 

If $\dim (Y)=0$ there is nothing to prove, since $\cI$ is trivial hence of order $0$. We assume proven order reduction in dimension $<\dim (Y)$.

We cover $V(\cI,a)$ with open patches $U$ possessing maximal contact hypersurfaces $H_U$. The coefficient ideal $C(\cI,a)_{|H_U}$ has order $\geq a!$. 

If this order is infinite, it means that $\cI_{|U}= \cI_{H_U}^a$, we simply blow up $H_U$ and automatically the order of $\cI$ is reduced on $U$. 

Otherwise we can inductively reduce  the order of this ideal by a functorial sequence of transformations $H_k \to \dots \to H$  until the order drops  below $a!$. By Corollary \ref{Cor:lifting} these provide a local order reduction for $(\cI,a)$ which patches together to a functorial order reduction by Proposition \ref{Prop:loca-aut}.

\subsection{The re-embedding principle}\label{Sec:reembedding} 
I still need to explain why a suitable embedded resolution of singularities implies resolution in general. If $X$ is covered by open subsets $X_i$ embedded in regular varieties $Y_i$, we want to claim that the resolutions $X'_i \to X_i$ agree on intersections $X_i\cap X_j$. Said another way, no matter how $X_i\cap X_j$ is embedded, the resolutions agree\footnote{There is a subtle issue of synchronization by codimension I will ignore. See \cite[\S 5.3]{BMT}, \cite[\S2.5.10]{Temkin-absolute}}. Since our procedures are functorial for \'etale maps and $Y_i$ are regular, we may as well assume $Y_i = \AA^{n_i}$. Finally affine spaces differ by iterated projections, so we are reduced to the following statement, which seems to follow from our procedures ``by accident'', see \cite[Claim 3.71.2]{Kollar}:

\begin{proposition}[The re-embedding principle]\label{Prop:reembedding}
Suppose $\cI$ is an ideal on a regular variety $Y$. Consider the embedding $Y \subset Y_1 := Y\times \AA^1$ sending $y\mapsto (y,0)$. Let $\cI_1 =  \cI \cO_{Y_1} + (z)$, where $z$ is the coordinate on $\AA^1$. Then the principalization described above of $\cI_1$ on $Y_1$ is obtained by taking the principalization of $\cI$ on $Y$ and blowing up the same centers, embedded in $Y_1$ and is transforms.
\end{proposition}

I would very much like to say that this follows from functoriality, but this is not so simple (see Koll\'ar's treatment). Instead, we look under the hood of principalization. To principalize $\cI_1$ we need to reduce the order of $\cI_1$ below 1. The order of $\cI_1$ is $1$, since $z$ has order 1, and then $z$ defines a maximal contact hypersurface, which is, seemingly by accident, precisely $Y$ with the coefficient ideal $\cI$, so the statement follows for the first blowing up. The local product structure persists after blowing up, so the statement  holds for the entire order reduction procedure.

\section{Toric varieties and toroidal embeddings} To proceed further it is useful to introduce a nice class of variety with ``fairly simple'' singularities.

\subsection{Toric varieties} A \emph{toric variety} is a \emph{normal} variety $X$ with a dense embedding $T = (\CC^*)^n \hookrightarrow X$ such that the action of $T$ on itself by translations extends to $X$. Here \emph{normal} means that the local rings are integrally closed, a condition which guarantees that $X$ is regular in codimension 1.  

Toric varieties are a simple playing ground for algebraic geometers, as many aspects of a toric variety can be translated to combinatorics. To a toric variety $X$ one associates a \emph{fan} $\Sigma_X$, a collection of rational polyhedral cones in the lattice $N_T:= \Hom(\CC^*, T)$ which intersect similarly to cells of  a CW complex. One makes toric varieties into a category on which arrows are torus equivariant morphisms which are surjective on the tori. Similarly fans form a category: a map of fans $\Sigma_1 \to \Sigma_2$ is induced by a map of lattices $N_{1} \to N_{2}$ with finite cokernel, such that a cone of $\Sigma_1$ maps into a cone in $\Sigma_2$.  There is an equivalence of categories 
$$\{\text{toric varieties}\} \leftrightarrow \{\text{fans}\}.$$

A toric variety $X$ is regular if and only if its fan $\Sigma_X$ is regular: every cone is simplicial, and generators of its edges span a saturated lattice in $N_T$. Toric birational maps correspond to subdivisions of fans, and so toric resolution of singularities can be done by finding a regular subdivision, a fairly simple task.

There are great sources to learn the theory. See \cite{KKMS, Oda, Fulton}.

\subsection{Toroidal embeddings} All toric varieties are rational, so they have a limited chance to help with resolution of singularities. A \emph{toroidal embedding} is an open embedding $U\subset X$ which locally analytically in the euclidean topology looks like a toric variety: for a point $p\in X$ there is a patch $V_p\subset X$ and a corresponding open set $W_p \subset Y$, where $T \subset Y$ is a toric variety, and an analytic isomorphism $V_p \to W_p$ carrying $V_p\cap U$ onto  $W_p \cap T$.

One can speak of \emph{toroidal morphisms} $X_1 \to X_2$ between toroidal embeddings: these are those morphisms which locally on the source look like toric morphisms of toric varieties.

As toroidal embeddings look locally like toric varieties their singularities are toric. It comes as no surprise that toric resolution of singularities extends quite easily to toroidal embeddings. In fact, one associates to a toroidal embedding $U \subset X$ a combinatorial gadget - a rational polyhedral cone complex - in a functorial manner. This is not an equivalence of categories, but it is still true that subdivisions correspond to toroidal birational morphisms, and the resolution procedure for fans extends to polyhedral cone complexes.

This is developed in \cite{KKMS}.


\section{Resolution in families}
I briefly recall known results on resolution in families, all relying on de Jong's \emph{alteration method.}
\subsection{The alteration theorem}
In \cite{dJ}, Johan de Jong discovered a method to replace a variety  $X$ by a regular variety $X'$ with a morphism $X' \to X$ which is not necessarily birational, but is \emph{proper, surjective and generically finite}. Such maps he called \emph{alterations}, which differ from birational \emph{modifications} in that the extension of function fields $K(X) \subset K(X')$ may be nontrivial.

\begin{theorem}[{\cite{dJ}}]\label{Th:alteration}
Let $X$ be a variety over a field of arbitrary characteristic and $Z$ a subvariety. There is an alteration $f:X' \to X$ such that $X'$ is smooth and $f^{-1}Z$ a simple normal crossings divisor.
\end{theorem}


\subsubsection{Sketch of proof} Here is the basic idea: assume for simplicity $X$ is projective; blowing up makes $Z$ into a divisor. One can choose a rational projection $X \das \PP^{n-1}$, which becomes a morphism after replacing $X$ by a modification, so that the generic fiber $X_\eta$ over the generic point $\eta \in \PP^n$ is a smooth curve, of some genus $g$, and $Z$ can be viewed as a collection of $k$ marked points on $X_\eta$. This corresponds to a morphism $\{\eta\} \to \ocM_{g,k}(X,d)$, the Kontsevich space of stable maps, where $d$ is the degree of $X_\eta$ with respect to some projective embedding. Properness of this moduli space provides us an alteration $B \to \PP^n$ over which this extends to a family of stable maps $Y_0 \to X$ parametrized by $B$. Induction on the dimension allows us to assume that $B$ is smooth and the degeneracy locus of $Y_0/B$ is a simple normal crossings divisor. An inspection of $Y_0$ shows that it has the structure of toroidal embedding, hence admits a combinatorial resolution of singularities $Y \to Y_0$. The composite morphism $Y \to X$ is the required alteration.

\subsection{Toroidalization} In the introduction, we stated resolution of singularities as the problem of making points of $X$ \emph{simple}. If instead we have a family of varieties $X\to B$ parametrized by a variety $B$, what should a resolution of singularities of the \emph{family} mean? That is, when are the singularities of the family \emph{simple}? 

It is not hard to see that making all the fibers regular is impossible. Since we agree that toroidal singularities are rather simple, one might consider a toroidal morphism to represent a family with simple singlarities. Here are two solutions based on this idea:

\begin{theorem}[Altered toroidalization, {\cite{dJ}}]\label{Th:alt-toroidalization} Let $X \to B$ be a dominant morphism of varieties. There are  alterations $B_1 \to B$ and $X_1 \to X\times_BB_1$, with regular toroidal embedding structures $U_B \subset B_1$ and $U_X \subset X_1$, and a \emph{toroidal} morphism $X_1 \to B_1$ making the following diagram commutative.
\begin{equation}\label{eq:alteration} \xymatrix{X_1 \ar[r]\ar[d] & X\ar[d] \\ B_1 \ar[r] & B}\end{equation}
\end{theorem} 

See improvement on this in Theorem \ref{Th:alt-semistable} below. This is proven in much the same way as one proves de Jong's alteration theorem, Theorem \ref{Th:alteration}. One needs to simply replace the projection $X \das \PP^n$ by a relative projection $X \das \PP^{d-1} \times B$, where $d$ is the relative dimension of $X$ over $B$.

In characteristic 0 one can inprove the situation, using \emph{modifications} instead of \emph{alterations}:

\begin{theorem}[Toroidalization, {\cite[Theorem 2.1]{AK}}]\label{Th:toroidalization} Let $X \to B$ be a dominant morphism of complex projective varieties. There is a modification $B' \to B$, a modification $X' \to X$, and regular toroidal embedding structures $U_B \subset B'$ and $U_X \subset X'$, such that the map $X' \das B'$ is a \emph{toroidal} morphism.
\end{theorem} 


Theorem \ref{Th:toroidalization} is proven using the following addition to de Jong's method, introduced in \cite{Abramovich-deJong}: in essence, one brings oneself to a situation  as in equation \eqref{eq:alteration}, where the Galois group $Gal(K(X_1)/ K(X))$ of the function field extension acts on the whole diagram. In characteristic 0 it turns out that the singularities of the quotients $X_1/G$ and $B_1/G$ can be resolved by toroidal methods  - this is a feature of \emph{tame} group actions in general. This sketch is only true in essense: in practice the Galois structure is intertwined with the inductive structure of the proof of   Theorem \ref{Th:alt-toroidalization}.

Theorems \ref{Th:alt-toroidalization} and \ref{Th:toroidalization}   have two major disadvantages: they are by no means functorial, and they necessarily change the general fiber of $X \to B$ even if it is already regular.

\subsection{Semistable reduction}
We have already pointed out that toroidal embeddings can be resolved. The same is true to some extent for families as well. This means that the singularities in Theorems \ref{Th:alt-toroidalization} and \ref{Th:toroidalization} can still be improved.

Let $X \to B$ be a toroidal morphism between toroidal embeddings. Assume $X$ and $B$ are regular. We say that $X \to B$ is \emph{semistable} if locally at any point $x\in X$ there are distinct monomial variables $y_1,\ldots,y_k$ on $B$ and $x_1,\ldots, x_m$ on $X$ such that $X$ is described by the following equations:

\begin{eqnarray*}x_1\cdots x_{l_1} &= y_1, \\ x_{l_1+1}\cdots x_{l_2} &= y_2,\\ \vdots \quad & \vdots \\ x_{l_{k-1}+1}\cdots x_{l_k} &= y_k.\end{eqnarray*}

This means that locally $X$ is a product of families of the form $x_1\cdots x_l = y$. This is truly the best one can hope for. A somewhat weaker and more flexible version would replace  $y_i$ by monomials $m_i$ without common factors. 

De Jong actually proved:

\begin{theorem}[Altered semistable reduction, {\cite{dJ}}]\label{Th:alt-semistable} Let $X \to B$ be a dominant morphism of varieties over a field. There are  alterations $B_1 \to B$ and $X_1 \to X\times_BB_1$, and regular toroidal embedding structures $U_B \subset B_1$ and $U_X \subset X_1$, with semistable morphism $X_1 \to B_1$. \end{theorem} 

In characteristic 0 a somewhat weaker result was proven in \cite[Theorem 0.3]{AK} where the generic fiber is modified but not altered. The strongest form is given in \cite{Karu-semistable} for families of surfaces and threefolds:

\begin{theorem}[Semistable reduction, {\cite{Karu-semistable}}]\label{Th:semistable} Let $X \to B$ be a dominant morphism of complex projective varieties with $\dim(X) - \dim(B) \leq 3$. There is an alteration $B' \to B$, a modification $X' \to X\times_B B'$, and regular toroidal embedding structures $U_B \subset B'$ and $U_X \subset X'$, such that the map $X' \das B'$ is semistable.
\end{theorem} 

The case of arbitrary relative dimension is conjectured in \cite[Conjecture 0.2]{AK}, and reduced to a completely combinatorial problem in \cite[Conjecture 8.4]{AK}. This remains open.

\section{Resolution in toroidal orbifolds}\label{Sec:ATW}

\subsection{Towards functorial resolution of families}

All the toroidalization and semistable reduction theorems above suffer from severe non-functoriality, and most importantly they change the generic fiber even if it is smooth. This is a major drawback in application. For example, one would like to take a smooth family of varieties over an open base and compactify it with as simple fibers as possible. One envisions compactifying by closure in some projective space and applying toroidalization or semistable reduction. The theorems above do not provide this, as the original family is necessarily changed.

The approach I present here is to start from scratch and use Hironaka's method instead of de Jong's. People have thought of this for a while, notably Cutkosky, see \cite{Cutkosky-m}, though his goals are different. 

\subsection{Temkin's functoriality principle} 
Consider the family $X \to B$ where $X$ is a regular surface with coordinates $x,y$ and $B$ a curve with coordinate $t$, and where the map is given by $xy = t$. This is a semistable family. 

Now take the base change $B_1 \to B$ given by $s^2 = t$. The pullback family $X_1 \to B_1$ is given by equation $xy = s^2$, which is semistable in the weaker sense, as $s^2$ is a monomial. If we are to allow semistable families to be compatible with base change this additional flexibility is a must. From the point of view of resolution of singularities in families, both families are good, even though $X_1$ is singular.

An important point in this example is that $B_1 \to B$ and $X_1 \to X$ are \emph{toroidal} morphisms. If we are to prove a functorial procedure for resolving singularities in families, the procedure must not modify families which are already semistable, so both $X \to B$ and $X_1 \to B_1$ must stay intact. 

In Hironaka's resolution, the best way to ensure that regular varieties stay intact is to require the resolution to be functorial  for smooth morphisms. Indeed if $X$ is regular then $X \to \Spec \CC$ is a smooth morphism, so the resolution of $X$ must be the pullback of the resolution of $\Spec \CC$, which is necessarily trivial.

In the semistable reduction problem, the best way to ensure that a toroidal $X \to B$ stays intact is to require functoriality for \emph{toroidal} morphisms. For instance, if already $B = \Spec \CC$ is just a point and $X$ is toroidal, then the morphism $X \to \Spec \CC$ is toroidal, so the procedure we produce for $X \to \Spec \CC$ must be the pullback of the procedure we produce for the identity $\Spec \CC \to \Spec \CC$, which is necessarily trivial. 

Temkin's functoriality principle is thus: 
\begin{quote} \em Toroidal morphisms form the smallest reasonable class of morphisms under which semistable reduction should be functorial. \end{quote}
And, in view of W\l odarczyk's philosophy,
\begin{quote} 
\em functoriality for toroidal morphisms should lead one to discover a semistable reduction procedure.
\end{quote}

As a necessary prerequisite, \emph{we must produce a resolution procedure, built on Hironaka's procedure, which is functorial for toroidal morphisms.}

\subsection{Enter toroidal orbifolds} 
Having accepted Temkin's functoriality principle and agreed that we must produce a resolution procedure, built on Hironaka's procedure, which is functorial for toroidal morphisms, there is another surprising conclusion coming our way.

Temkin's principle forces us to take a departure from previous algorithms:

\begin{enumerate}
\item We can no longer work with a smooth ambient variety $Y$ - \emph{we must allow $Y$ to have toroidal singularities.}
\item We cannot use only blowings up of smooth centers as our basic operations. Instead we use a class of modifications, called \emph{Kummer blowings up}, which are stable under toroidal base change. These involve taking roots of monomials, in particular:
\item We can no longer work only with varieties $Y$ - \emph{we must allow $Y$ to be a Deligne--Mumford stack.}
\end{enumerate}

In essence, we are using ``weighted blowings up on steroids''.  The stacks we need are as follows:

\begin{definition}
A \emph{toroidal orbifold} $(Y,U)$ is a Deligne--mumford stack $Y$ with diagonalizable inertia with a toroidal embedding $U\subset Y$.
\end{definition}

Item (3) may be hard to accept but it is absolutely essential for functoriality under toroidal morphisms. 

Consider the affine plane $Y$, with coordinates $x,u$, where we endow the plane with a toroidal structure by declaring $U = Y \setminus \{u=0\}$, so $x$ is a parameter and $u$ is a \emph{monomial}. Say we want to principalize the ideal $\cI:=(x^2, u^2)$. Functoriality under smooth morphisms suggests that we must blow up the ideal $(x,u)$. This indeed works and principalizes  $\cI$ in one step.

Now consider the affine plane $Y_0$, with coordinates $x,v$, similar to the above, but say we want to principalize the ideal $\cI_0:=(x^2,v)$. Note that we have a toroidal morphism $Y \to Y_0$ given by $u^2 = v$, and $\cI = \cI_0 \cO_Y$. Temkin's functoriality tells us that there should be a center on $Y_0$ whose pullback is $(x,u)$. This center must therefore be defined by $(x,\sqrt v)$!

There is only one way to deal with it, and that is to work systematically in a setup where one is allowed, when necessary, to take roots of monomials. This is possible precisely when working with toroidal orbifolds, as indicated above.

\subsection{Principalization in toroidal orbifolds} In the joint work \cite{ATW-principalization} with Temkin and W\l odarczyk we prove
\begin{theorem}[Toroidal principalization]\label{Th:principalization}
Let $(Y, U)$ be a toroidal orbifold and $\cI$ an ideal sheaf. There is a sequence  $Y_n \to Y_{n-1} \to \dots\to Y_0 = Y$ of  Kummer blowings up, all supported over the vanishing locus $V(\cI)$, such that $\cI \cO_{Y_n}$ is an invertible monomial ideal. \emph{The process is functorial for toroidal  base change morphisms  $Y'\to Y$.}
\end{theorem}

\subsection{Logarithmic derivatives and logarithmic orders}
Temkin's functoriality principle suggests a natural replacement for derivatives.

In the toroidal world, the natural replacement for derivatives is provided by \emph{logarithmic derivatives}: if $u$ is a \emph{monomial} function on a toroidal $Y$, then we use the operator $u \frac{\partial}{\partial u}$, which sends $u$ to itself, but not $\frac{\partial}{\partial u}$. One then defines the sheaf $\cD^{\leq a}_{Y,U}$ of \emph{logarithmic} differential operators of order $\leq a$. 

Given an ideal $\cI$ on a toroidal $(Y,U)$, one defines its \emph{logarithmic order} at $p\in Y$ to be 
$\logord_{p}(\cI) = \min\{a: \cD^{\leq a}_{Y,U}(\cI) = (1)\}.$ This can take the value $\infty$ when the \emph{monomial part} $\cM(\cI):= \cD^{\infty}_{Y,U}\cI$ is nontrivial.

Then a miracle happens: using Temkin's functoriality and logarithmic derivatives, the broad outlines of principalization described above, as laid out in detail in  \cite{Wlodarczyk}, work in this new context 
%
%
 once one has a stable formalism of toroidal orbifolds. The formalism is developed in \cite{ATW-destackification} and the proof of the theorem  is written out in \cite{ATW-principalization}.

Our next task is to return to work on families of varieties. We hope to report on that in the near future.

\bibliographystyle{amsalpha}
\bibliography{ICM}

\end{document}